\documentclass[10pt]{article}
\usepackage{amssymb,amsmath,textcomp}
\usepackage{color}
\newtheorem{defi}{Definition}[section]
\newtheorem{satz}{Theorem}[section]
\newtheorem{theorem}{Theorem}[section]
\newtheorem{lemma}{Lemma}[section]

\newtheorem{remark}{Remark}[section]
\newtheorem{bei}{Example}[section]
\newtheorem{korollar}{Corollary}[section]

\newcommand{\iR}{\mathbb{R}}
\newcommand{\iN}{\mathbb{N}}

\newcommand{\iC}{\mathbb{C}}

\newcommand{\dif}{\mathrm{d}}

\begin{document}
\begin{center}
{\bf\Large Long-time behaviour of non-local in time Fokker-Planck equations via the entropy method}
\end{center}
\vspace{0.7em}
\begin{center}
Jukka Kemppainen\footnote{J.K. was supported by the Vilho, Yrj\"{o} and Kalle V\"{a}is\"{a}l\"{a} Foundation of the Finnish Academy of Science and Letters} and Rico Zacher\footnote{R.Z.\ was supported by a research grant of the German Research Foundation (DFG), GZ  Za 547/4-1.}
\end{center}
\begin{abstract}
We consider a rather general class of non-local in time Fokker-Planck equations and show by means
of the entropy method that
as $t\to \infty$ the solution converges in $L^1$ to the unique steady state. 
Important special cases are the time-fractional and ultraslow diffusion case.
We also prove estimates for the rate of decay. In contrast to the classical (local) case, where
the usual time derivative appears in the Fokker-Planck equation, the obtained decay rate depends on the entropy, which is related
to the integrability of the initial datum. It seems that higher integrability of the initial
datum leads to better decay rates and that the {\em optimal decay rate} is reached, as we show, when 
the initial datum belongs to a certain weighted \(L^2\) space. We also show how our estimates can be
adapted to the discrete-time case thereby improving known decay rates from the literature. 
\end{abstract}
\vspace{0.7em}
\begin{center}
\end{center}

\section{Introduction and main result}
We study the large time behaviour of solutions to the non-local in time Fokker-Planck
equation
\begin{equation} \label{FP_eq}
\partial_t \big(k\ast [u-u_0]\big)-\nabla\cdot\left(\nabla u+u\nabla V\right)=0,\quad
t>0,\,x\in \iR^d,
\end{equation}
where \(u=u(t,x)\),
\begin{equation} \label{FP_init_cond}
u|_{t=0}=u_0,\,x\in \iR^d,
\end{equation}
is the initial datum, \(V=V(x)\) is the given potential and $k\ast v$ denotes the convolution on the positive halfline $\iR_+:=[0,\infty)$ w.r.t.\ the time variable,
that is
$(k\ast v)(t)=\int_0^t k(t-\tau)v(\tau)\,d\tau$, $t\ge 0$. 

The kernel $k\in L^1_{loc}(\iR_+)$ belongs to a certain class of kernels (it is of type ${\cal PC}$, see 
Section \ref{SecKernels}), which covers most of the
relevant integro-differential operators w.r.t.\ time that appear in physics in the context of so-called {\em subdiffusion}
processes, in particular the {\em time-fractional} case, where 
\begin{equation} \label{kexample}
k(t)=\frac{t^{-\alpha}}{\Gamma(1-\alpha)}, \quad t>0,
\end{equation} and $\partial_t(k\ast v)$ coincides with 
the Riemann-Liouville fractional derivative $\partial_t^\alpha v$ of order $\alpha\in (0,1)$, see 
Section \ref{SecKernels} below. Note that $\partial_t(k\ast[v-v(0)])=k\ast \partial_t v$
for sufficiently smooth functions $v$.

Fokker-Planck equations, also known as {\em Kolmogorov forward equations}, are a central subject
in the theory of diffusion processes. They describe how the probability density function of 
the velocity (or some other observable like position) of a particle evolves with time.
The non-local in time version~\eqref{FP_eq} is used in physics for the modeling of subdiffusion processes.
Here the mean square displacement $m(t)$ of a diffusive particle, which is defined by
\begin{equation} \label{MSDdef}
m(t)=\int_{\iR^d}|x|^2 Z(t,x)\,\dif x,\quad t>0,
\end{equation}
with $Z(t,x)$ denoting the fundamental solution of equation~\eqref{FP_eq} with $V=0$ and 
$Z|_{t=0}=\delta_0$ (see \cite{KempSiljVergZach16}), grows slower as $t\to \infty$ than in the classical case of Brownian motion, where $m(t)=ct$, $t>0$,
with some constant $c>0$. For example, in the time-fractional diffusion case $m(t)=ct^\alpha$ (cf.\ \cite{Metz}). Another example, which is also covered by our setting, is ultraslow diffusion (see Example \ref{beislow} below);
here $m(t)$ merely grows logarithmically for large times. 

In the time-fractional case, equation~\eqref{FP_eq} may be viewed as the master equation for the probability density of the scaling limit of a random walk with a heavy-tailed waiting time distribution
or as an overdamped limit of a fractional Klein-Kramers equation. For details we refer to~\cite{Metz,Metz2}. 
Concerning the stochastic viewpoint of~\eqref{FP_eq} in this case, we first recall that 
the stochastic differential equation (SDE) corresponding to the classical Fokker-Planck equation
\begin{equation} \label{FP_classic}
\partial_t u-\nabla\cdot\left(\nabla u+u\nabla V\right)=0,\quad
t>0,\,x\in \iR^d,
\end{equation}
is given by
\begin{equation}\label{SDE_Fokker_Planck}
\dif X_t=-\nabla V(X_t)\dif t+\sqrt{2}\dif B_t, 
\end{equation}
where \(B_t\) denotes Brownian motion. Let \(D_t\) be an \(\alpha\)-stable scaling limit of
the waiting time process and let
\[
E_t=\mathrm{inf}\{\tau>0\,:\,D_\tau>t\}
\]
be its inverse or the first hitting time process. Then the time-changed
process \(Y_t=X_{E_t}\) satisfies the SDE
\begin{equation}\label{SDE_nonlocal_Fokker_Planck}
\dif Y_t=-\nabla V(Y_t)\dif E_t+\sqrt{2}\dif B_{E_t} 
\end{equation}
and the transition probability density \(u\) is the solution of~\eqref{FP_eq} with 
$k$ as in \eqref{kexample}, see~\cite{HKU,MNV} for more details.

The slowness of subdiffusion processes for large times is reflected by lower temporal decay rates as well.
These have been studied extensively in~\cite{KempSiljVergZach16,KempSiljZach17} for solutions of ~\eqref{FP_eq} with \(V\equiv 0\). For example, considering the time-fractional case (that is $k$ is given by \eqref{kexample}
with $\alpha\in (0,1)$) and assuming that $u_0\in L^1(\iR^d)\cap L^2(\iR^d)$ the solution $u$ of~\eqref{FP_eq} with vanishing potential given by 
\[
u(t,x)=\int_{\iR^d} Z(t,x-y)u_0(y)\,\dif y
\]
satisfies the (optimal) estimates 
\[
||u(t,\cdot)||_{L^{2}(\iR^d)}\lesssim t^{-\min\{\frac{\alpha d}{4},\alpha\}},\quad t>0,\;d\in \iN\setminus \{4\}
\]
and
\[
||u(t,\cdot)||_{L^{2,\infty}(\iR^d)}\lesssim t^{-\alpha},\quad t>0,\;d=4,
\]
where $L^{2,\infty}(\iR^d)$ refers to the weak $L^2$ space, see \cite{KempSiljVergZach16}. This 
decay behaviour is markedly different from that in the case of the classical heat equation, where
under the same assumptions $||u(t,\cdot)||_{L^{2}(\iR^d)}\lesssim t^{-d/4}$, $t>0$, in all dimensions $d\in \iN$.
 
If the potential $V$ in~\eqref{FP_eq} does not vanish and $e^{-V}\in L^1(\iR^d)$ (which will be 
assumed throughout this paper), there is a unique steady state of~\eqref{FP_eq} with unit mass,
which is given by
\begin{equation}\label{Gibbs_distribution}
u_\infty(x)=M\mathrm{e}^{-V(x)},\quad M=\left(\int_{\iR^d}\mathrm{e}^{-V(y)}\dif y\right)^{-1}.
\end{equation}
A typical choice of the potential is 
\[
V(x)=\frac{1}{2}m|x|^2,\quad m>0,
\]
which corresponds to kinetic energy when regarding $x$ as the velocity of the particle. In this case
the steady state becomes a Gaussian. 

Assuming that $u_0$ is a probability density we aim at showing that the solution of 
the initial-value problem \eqref{FP_eq}, \eqref{FP_init_cond} tends to $u_\infty$ in $L^1(\iR^d)$
as $t\to \infty$, and we are interested in convergence rates which are as precise as possible.
We should point out that we do not address existence and regularity of solutions to \eqref{FP_eq}
but perform formal calculations, which become rigorous provided the solution is sufficiently smooth 
(and satisfies a suitable growth condition for $|x|\to \infty$) or
can be justified by appropriate regularization techniques.

As to existence and regularity of solutions to \eqref{FP_eq}, we are not aware of any work where this specific
problem has been studied. However, there is a theory of abstract evolutionary integral equations, see
e.g.\ \cite{JanI,ZEQ,ZWH}, which generalizes semigroup theory and can be applied to \eqref{FP_eq}. For example,
by the {\em subordination principle} from Pr\"uss \cite[Chapter 4]{JanI}, any result on the generation of a $C_0$-semigroup for the classical Fokker-Planck equation \eqref{FP_classic} yields
a related well-posedness result for equation \eqref{FP_eq} if the kernel $k$ is of type ${\cal PC}$. The theory in \cite{ZWH}
can be used to prove existence and uniqueness of a weak solution to \eqref{FP_eq} in the Hilbert space setting considered
in Section \ref{optHilbert} below. Higher regularity of the solution can be obtained by means of the results and techniques from
\cite[Chapter 2,3]{JanI} (maximal regularity in H\"older spaces in time, classical solutions) and \cite{ZEQ} (strong $L^p$-solutions), at least in the time-fractional case. An analogue of the celebrated theorem of Jordan, Kinderlehrer and Otto
on the gradient flow structure of the classical Fokker-Planck equation in the Wasserstein space ${\cal P}_2(\iR^d)$ (\cite{JKO})
seems to be unknown for equation \eqref{FP_eq} and would be highly desirable.

For the classical Fokker-Planck equation \eqref{FP_classic}
and certain variants of it (including nonlinear problems), there is a huge literature on
convergence to equilibrium as $t\to \infty$, see e.g.\ \cite{AMTU,jungel2016,MV,To} and the references given therein. A very powerful technique in this context is the entropy method as presented, e.g., in~\cite{jungel2016}. The key idea
consists in proving a differential inequality of the form
\begin{equation}\label{diff_ineq_entropy}
\frac{\dif}{\dif t}H(u(t))\le -\kappa H(u(t)),\quad t>0,
\end{equation}
with $\kappa>0$, for a suitable relative entropy \(H\) associated with the steady state $u_\infty$.
Such a relative entropy is of the form
\begin{equation}\label{relative_entropy}
H(u):=H_\phi(u)=\int_{\iR^d}\phi\left(\frac{u(x)}{u_\infty(x)}\right)u_\infty(x)\dif x,\quad 
u\in L^1(\iR^d),\;u\ge 0,\;||u||_1=1,
\end{equation}
with {\em generating function} $\phi:\,[0,\infty)\rightarrow \iR$, which is strictly convex and satisfies
$\phi(1)=0$. Note that Jensen's inequality implies that $H(u)\ge 0$. An important example is given by
\begin{equation}\label{entropy_function_p}
\phi_\beta(x)=x^\beta-1-\beta(x-1),\quad 1<\beta\le 2.
\end{equation}
Inequality~\eqref{diff_ineq_entropy}
then implies the exponential decay of the relative entropy,
\[
H(u(t))\le H(u(0))\mathrm{e}^{-\kappa t},\quad t>0,
\]
which together with a suitable Csiz\'ar-Kullback-Pinsker inequality 
(see Theorem \ref{CKP_ineq} below) leads to an
exponential decay estimate for $||u(t)-u_\infty||_{L^1(\iR^d)}$ provided that $H(u(0))<\infty$. It is well known (see e.g.\
\cite{AMTU} and \cite[Chapter 2]{jungel2016}) that \eqref{diff_ineq_entropy} with $\kappa=2\lambda$ can be established for smooth potentials satisfying the convexity
condition (also termed {\em Bakry-Emery condition})
\begin{equation}\label{potential_convexity_cond}
\nabla^2 V(x)\ge \lambda I,\quad x\in \iR^d,
\end{equation}
where $\lambda>0$ and $\nabla^2 V$ denotes the Hessian of $V$. This is possible since 
condition \eqref{potential_convexity_cond} ensures the validity of certain convex Sobolev inequalities,
in particular the logarithmic Sobolev inequality and a weighted Poincar\'e inequality.

Assuming that condition \eqref{potential_convexity_cond} is satisfied,
the key idea of the entropy method in the non-local in time case is to derive an integro-differential inequality
of the form
\begin{equation}\label{integro_diff_ineq}
\frac{\dif}{\dif t}(k*[\Phi-\Phi_0])(t)+\mu \Phi(t)\le 0,\quad t>0,
\end{equation}
where $\mu>0$, $\Phi_0=\Phi(0)$ and $\Phi(t)$ is a suitable function depending on a relative entropy, e.g.\ a power $\Phi(t)
=H(u(t))^{1/\beta}$ with $\beta\in [1,2]$. By the comparison principle 
(cf.~\cite[Section 2.3]{VZ}), \eqref{integro_diff_ineq} then implies the estimate
\begin{equation}\label{integro_diff_ineq_corol}
\Phi(t)\le s_\mu(t)\Phi_0,
\end{equation}
where \(s_\mu\), the so-called {\em relaxation function}, solves the problem
\begin{equation} \label{relaxfdef}
\frac{\dif}{\dif t}(k*[s_\mu-1])(t)+\mu s_\mu(t)\le 0,\quad t>0,\quad s_\mu(0)=1.
\end{equation}

The large-time behaviour of the relaxation function depends heavily on the kernel $k$. Denoting
by $\hat{k}$ the Laplace transform of $k$, it is known that $s_\mu(t)\to 0$ as $t\to\infty$ if and only
if 
\begin{equation} \label{kLTcond}
\lim_{\lambda\to 0+} \frac{1}{\lambda \hat{k}(\lambda)}=\infty.
\end{equation}  
This condition is satisfied for most of the important examples, in particular in the time-fractional case
\eqref{kexample} for all $\alpha\in (0,1)$; here $\hat{k}(\lambda)=\lambda^{\alpha-1}$, Re$\,\lambda>0$.

Our main result is the following. Here, by saying that the kernel $k$ is of type ${\cal PC}$ we mean that $k\in L^1_{loc}(\iR_+)$ is nonnegative and nonincreasing, and that there exists a kernel 
$l\in L^1_{loc}(\iR_+)$ such that
$k\ast l=1$ in $(0,\infty)$, cf. Section \ref{SecKernels}.


\begin{satz}\label{thm_main_result}
Let $k$ be a kernel of type ${\cal PC}$ and $V\in C^2(\iR^d)$ satisfying condition 
\eqref{potential_convexity_cond} for some $\lambda>0$ and such that $e^{-V}\in L^1(\iR^d)$.
Let further $u_0\in L^1(\iR^d)$ be a probability density. Defining $u_\infty$ as in
\eqref{Gibbs_distribution} the following statements hold true.

{\bf Part A (general entropy)} Let $\phi$ be an admissible generating function according to Definition
\ref{rel_entropy_def}
and let $H(u)$ be the corresponding relative entropy associated with $u_\infty$ given by
\eqref{relative_entropy}. Assume that $H(u_0)<\infty$ and that the positive solution $u$ of  
\eqref{FP_eq}, \eqref{FP_init_cond} is sufficiently smooth. Then
\begin{equation} \label{decay_rate_entropyA}
H(u(t))\le s_{2\lambda}(t) H(u_0),\quad t>0.
\end{equation}
Moreover,
\begin{equation} \label{decay_rate_solutionA}
||u(t)-u_\infty||_{L^1(\iR^d)}\le C_\phi \sqrt{s_{2\lambda}(t)} \sqrt{ H(u_0)},\quad t>0,
\end{equation}
with $C_\phi=\sqrt{2/\phi''(1)}$. In particular, if in addition $k$ satisfies \eqref{kLTcond}, then $u(t)\to u_\infty$ in $L^1(\iR^d)$ as $t\to \infty$.

{\bf Part B (power type entropy)} Let $\beta\in (1,2]$ and $\phi_\beta$ be defined as in 
\eqref{entropy_function_p}. Let $H(u)$ be the corresponding relative entropy associated with $u_\infty$.  Assume that $H(u_0)<\infty$ and that the positive solution $u$ of  
\eqref{FP_eq}, \eqref{FP_init_cond} is sufficiently smooth. Then
\begin{equation} \label{decay_rate_entropyB}
H(u(t))\le \left(s_{2\lambda/\beta}(t)\right)^{\beta} H(u_0),\quad t>0.
\end{equation}
Moreover,
\begin{equation} \label{decay_rate_solutionB}
||u(t)-u_\infty||_{L^1(\iR^d)}\le \sqrt{\frac{2}{\beta(\beta-1)}}\left(s_{2\lambda/\beta}(t)\right)^{\beta/2}\sqrt{H(u_0)},\quad t>0.
\end{equation}
\end{satz}
Part A of Theorem \ref{thm_main_result} says that, for any admissible entropy, \eqref{integro_diff_ineq_corol} is true with $\Phi(t)=H(u(t))$ and $\mu=2\lambda$. That is,
under the additional assumption \eqref{kLTcond},
the entropy decays at least as fast as the relaxation function $s_{2\lambda}(t)$, and 
$||u(t)-u_\infty||_{L^1(\iR^d)}$ is controlled by the square root of $s_{2\lambda}(t)$. Note that Part A covers the important {\em Boltzmann entropy} defined by~\eqref{Boltzmann_entropy} in Section~\ref{sec:addmissible_entropies}.

These results are the analogues of those for the classical Fokker-Planck equation
\eqref{FP_classic}, where the corresponding relaxation function solves the ODE
\[
\frac{\dif}{\dif t}s_\mu(t)+\mu s_\mu(t)=0,\; t>0,\quad s_\mu(0)=1,
\]
with $\mu=2\lambda$, that is $s_\mu(t)=e^{-\mu t}=e^{-2\lambda t}$. No matter which entropy
is used, the norm $||u(t)-u_\infty||_{L^1(\iR^d)}$ decays in the classical case (at least as fast) as
$e^{-\lambda t}\,(=\sqrt{s_{2\lambda} (t)})$ (up to some constant).

As Part B of Theorem \ref{thm_main_result} shows, the latter is no longer the case in the non-local
in time situation, in the sense that the decay rates improve for higher values of $\beta\in (1,2]$ when
using the entropy generated by the power type function $\phi_\beta$. Observe that for $\beta=2$
we even reach the best possible decay estimate for $||u(t)-u_\infty||_{L^1(\iR^d)}$ one can hope for, with
a factor $s_\lambda(t)$ on the right-hand side of \eqref{decay_rate_solutionB}. In the classical
case, we have the identity
\[
\left(s_{2\lambda/\beta}(t)\right)^{\beta/2}=\left(e^{-\frac{2\lambda t}
{\beta}}\right)^{\beta/2}=e^{-\lambda t}=s_\lambda(t),
\]
which is not true any more in the non-local case, where, in general,  $s_\lambda(t)$ decays faster than
$\left(s_{2\lambda/\beta}(t)\right)^{\beta/2}$ for $\beta\in [1,2)$.

Notice as well that
in case of the generating function $\phi_\beta$, the condition $H(u_0)<\infty$ is equivalent to
\[
\int_{\iR^d} \left(\frac{u_0}{u_\infty}\right)^\beta u_\infty\,\dif x<\infty,
\] 
which is stronger for larger values of $\beta$, by H\"older's inequality. Thus in contrast to
the local Fokker-Planck equation, the decay rates improve with higher integrability of the 
initial datum w.r.t. the corresponding weighted Lebesgue space.

In the time-fractional case, the decay behaviour of the relaxation function can be quantified 
precisely. Theorem \ref{thm_main_result} leads to the following result.
\begin{korollar} \label{main_result_timefrac}
{\bf (fractional dynamics)} Let $\alpha\in(0,1)$ and $k(t)=\frac{t^{-\alpha}}{\Gamma(1-\alpha)}$, $t>0$. Let $V,\lambda$, and
$u_\infty$ be as in Theorem \ref{thm_main_result} and suppose that
$u_0\in L^1(\iR^d)$ is a probability density.

{\bf Part A (general entropy)} Let $\phi$ be an admissible generating function
and $H(u)$ be the corresponding relative entropy associated with $u_\infty$. Assume that $H(u_0)<\infty$ and that the positive solution $u$ of  
\eqref{FP_eq}, \eqref{FP_init_cond} is sufficiently smooth. Then
\begin{equation} \label{decay_rate_entropyA_timefrac}
H(u(t))\le \frac{C(\alpha,\lambda)}{1+t^\alpha}\, H(u_0),\quad t>0,
\end{equation}
and
\begin{equation} \label{decay_rate_solutionA_timefrac}
||u(t)-u_\infty||_{L^1(\iR^d)}\le \frac{C(\alpha,\lambda,\phi)}{1+t^{\alpha/2}}\,
 \sqrt{ H(u_0)},\quad t>0.
\end{equation}

{\bf Part B (power type entropy)} Let $\beta\in (1,2]$ and $\phi_\beta$ be defined as in 
\eqref{entropy_function_p}. Let $H(u)$ be the corresponding relative entropy. Assume that $H(u_0)<\infty$ and that the positive solution $u$ of  
\eqref{FP_eq}, \eqref{FP_init_cond} is sufficiently smooth. Then
\begin{equation} \label{decay_rate_entropyB_timefrac}
H(u(t))\le \frac{C(\alpha,\beta,\lambda)}{1+t^{\alpha\beta}}\, H(u_0),\quad t>0,
\end{equation}
and
\begin{equation} \label{decay_rate_solutionB_timefrac}
||u(t)-u_\infty||_{L^1(\iR^d)}\le \frac{\tilde{C}(\alpha,\beta,\lambda)}{1+t^{\alpha\beta/2}}\,
\sqrt{H(u_0)},\quad t>0.
\end{equation}
\end{korollar}

Our proof of Theorem \ref{thm_main_result} makes heavily use of the fundamental
identity \eqref{fundamental_identity} for the integro-differential operator $Bv=\partial_t(k\ast v)$.
In order to obtain the improved estimates in the case of the power type entropy, we derive delicate estimates for the generating function $\phi_\beta$ and the associated entropy and employ
the fundamental identity in its full strength, that is, all the terms on the right hand side of~\eqref{fundamental_identity} are used. It is also noteworthy that we have to apply the fundamental identity twice.

It turns out that our method also applies to the {\em time-discretized} classical Fokker-Planck
equation. Using our estimates for the generating function $\phi_\beta$ and the associated entropy
we are seemingly able to improve existing results from the literature. The point is that in the time discrete
case (where the equation is non-local as well), the decay rates become better with higher
integrability of the initial datum, exactly as in Part B of Theorem \ref{thm_main_result}.
We strongly believe that our techniques may be also useful in the context of other time-discrete schemes
for classical PDEs.  

The paper is organized as follows. In Section 2 we collect some preliminaries on 
kernels of type ${\cal PC}$ and the associated relaxation functions, we recall the fundamental identity 
for operators of the form $\partial_t(k\ast \cdot)$, and we describe the class of admissible relative entropies and
recall the general Csisz\'ar-Kullback-Pinsker inequality as well as the convex Sobolev inequality.
Section 3 contains the proof of our main result, Theorem \ref{thm_main_result}. In Section 4 we discuss several
important examples of kernels of type ${\cal PC}$, in particular we consider the time-fractional case and prove
Corollary \ref{main_result_timefrac}. In Section 5 we discuss the optimality of the entropy decay rate in the 
case $\beta=2$. Finally, Section 6 is devoted to the time-discrete case.

\section{Preliminaries}

\subsection{Kernels of type ${\cal PC}$ and relaxation functions} \label{SecKernels}

Throughout this paper we assume that the kernel $k$ in \eqref{FP_eq} is of type ${\cal PC}$. 
By this we mean
that it satisfies the condition
\begin{itemize}
\item [{\bf ($\mathcal{PC}$)}] $k\in L^1_{loc}(\iR_+)$ is nonnegative and nonincreasing, and there exists a kernel 
$l\in L^1_{loc}(\iR_+)$ such that
$k\ast l=1$ in $(0,\infty)$.
\end{itemize}
In this case we also use the notation $(k,l)\in \mathcal{PC}$. Note that $(k,l)\in {\cal PC}$ implies that $l$ is nonnegative, cf.\ \cite[Theorem 2.2]{CN}. 

Condition ($\mathcal{PC}$) has already been used before in a couple of papers 
(see e.g.~\cite{KempSiljVergZach16,VZ,Za}). It covers many important examples of integro-differential operators w.r.t.\ time that
are used in physics to describe subdiffusion processes. An important example is given by $(k,l)=(g_{1-\alpha},g_\alpha)$ with $\alpha\in(0,1)$, where $g_\beta$ denotes the standard kernel
\begin{equation}\label{standard_kernel}
g_\beta(t)=\,\frac{t^{\beta-1}}{\Gamma(\beta)}\,,\quad
t>0,\quad\beta>0.
\end{equation}
In this case, the term $\partial_t(k\ast v)$ becomes the Riemann-Liouville fractional derivative
$\partial_t^\alpha v$, and $k \ast \partial_t v={}^c D_t^\alpha v$, the Caputo fractional derivative of order $\alpha$ (cf.\ \cite{KST}). We then also call (\ref{FP_eq}) {\em time-fractional Fokker-Planck equation}.

Condition ($\mathcal{PC}$) also contains the multi-term fractional diffusion case, 
see Example \ref{multiterm} below. 
Another interesting and important class of examples is given by
\begin{equation*}
k(t)=\int_0^1 g_\beta(t)\omega(\beta)\,\dif \beta,
\end{equation*}
where $\omega\in C([0,1])$ is a nonnegative weight function that does not vanish everywhere.
In this situation the operator $\partial_t(k\ast \cdot)$ is a so-called operator of
{\em distributed order}, and (\ref{FP_eq}) is an example of an {\em ultraslow Fokker-Planck equation} if $\omega(0)
\neq 0$ (\cite{Koch08}). As we shall see, the obtained decay rate for the entropy will partly justify the chosen name.
The special case $\omega\equiv 1$ is discussed in Example \ref{beislow} below.

Assuming $(k,l)\in \mathcal{PC}$ and $\mu\ge 0$ we recall that the associated relaxation function $s_\mu$ is defined 
via the problem \eqref{relaxfdef}. Convolving the first equation in \eqref{relaxfdef} with $l$, it is easy to check that
problem \eqref{relaxfdef} is equivalent to the Volterra equation
\[
s_\mu(t)+\mu (l\ast s_\mu)(t)=1,\quad t\ge 0.
\]
It is known that $s_\mu$ is nonnegative, nonincreasing, and that $s_\mu\in H^{1,1}_{loc}(\iR_+)$; furthermore
$\partial_\mu s_\mu(t)\le 0$, see e.g.\ Pr\"uss \cite{JanI}. Moreover, one knows that for any $\mu\ge 0$ there holds
\begin{equation} \label{smubounds}
\frac{1}{1+\mu \,k(t)^{-1}}\,\le s_\mu(t)\le \,\frac{1}{1+\mu \,(1\ast l)(t)}\,,\quad
t>0,
\end{equation}
which also entails that
\[
\big[1-s_\mu(t)\big] k(t)\,\le \mu s_\mu(t)\le \,\big[1-s_\mu(t)\big]\,\frac{1}{(1\ast l)(t)} \,,\quad t>0,
\]
see \cite{VZ}. This implies that for any fixed $\mu>0$, $s_\mu(t)$ cannot decay faster than the kernel $k(t)$, and $s_\mu(t)$ decays
at least like $(1\ast l)(t)^{-1}$. Note that $\lim_{t\to \infty} s_\mu(t)=0$ if and only if $l\notin L^{1}(\iR_+)$, see e.g.\
\cite[Lemma 6.1]{VZ}. Using
\[
\hat{k}(\lambda) \hat{l}(\lambda)= \frac{1}{\lambda},\quad \lambda>0,
\]
it is easy to see that $l\notin L^{1}(\iR_+)$ if and only if
condition \eqref{kLTcond} is satisfied.

\subsection{The fundamental identity for the operator $\partial_t(k\ast \cdot)$}

To simplify the notations, in what follows we use the symbol \(\partial_t\) for the derivative also for the functions that depend only on \(t\).
An important tool in our approach is the so-called
{\em fundamental identity} for integro-differential
operators of the form $\partial_t(k\ast \cdot)$, cf.\ also
\cite{Za,Za1}.
Let $k\in L^1_{loc}(\iR_+)$ be a nonnegative and nonincreasing kernel, $I\subset \iR$ an interval, $\psi \in C^1(I)$, 
and $u\in L^1((0,T))$ with $u(t)\in I$
for a.a. $t\in (0,T)$. Then, for sufficiently smooth $u$ and $k$ there holds for a.a.\ $t\in (0,T)$
\begin{equation}\label{fundamental_identity}
\begin{split}
\psi'(u(t))\partial_t (k*u)(t) &=\partial_t (k*\psi(u))(t)+\Big(-\psi(u(t))+\psi'(u(t))u(t)\Big)k(t)\\
&+\int_0^t \Big(\psi(u(t-s))-\psi(u(t))-\psi'(u(t))[u(t-s)-u(t)]\Big)(-k'(s))\dif s.
\end{split}
\end{equation}
This can be shown by a straightforward computation. Observe that the third term on the right-hand side is nonnegative if $\psi$ is convex. Assuming in addition that $\psi$ is convex and that $u_0\in I$, it follows from \eqref{fundamental_identity} that
\begin{equation} \label{FIconvex}
\psi'(u(t))\partial_t \big(k*[u-u_0]\big)(t)\ge \partial_t \big(k*[\psi(u)-\psi(u_0)]\big)(t),\quad t\in (0,T),
\end{equation}
cf.\ \cite[Corollary 6.1]{KempSiljVergZach16}. Note that here it is not assumed that $u_0=u(0)$.

We point out again that in this paper we perform formal calculations, that is, we assume that the solution
is smooth enough so that \eqref{fundamental_identity} and \eqref{FIconvex} can be used with a kernel
$k$ of type ${\cal PC}$, which is always singular at $t=0$. For less regular solutions (e.g.\ weak solutions in 
a certain sense), an important regularization technique consists in replacing the operator 
$\partial_t(k\ast \cdot)$ by its Yosida approximations $\partial_t(k_n\ast \cdot)$, $n\in \iN$, with the
more regular kernels $k_n=n s_n\in H^{1,1}_{loc}(\iR_+)$. We refer e.g.\ to \cite{VZ, Za,Za1}, where this
method has been used to prove rigorous estimates in a weak setting.

The following technical lemma is an application of the convexity inequality \eqref{FIconvex}.
The function $w$ and the kernel $k$ involved have to be sufficiently regular so that \eqref{FIconvex} can be used.
\begin{lemma} \label{dividingprefactor}
Let $T>0$ and $k\in L^1_{loc}(\iR_+)$ be a nonnegative and nonincreasing kernel. Let $w:\,[0,T]\rightarrow [0,\infty)$
be a sufficiently smooth function, $w_0\in [0,\infty)$ and $\delta, \gamma,\eta>0$ such that $\delta<\eta$. Let further
$\mu\in L^1((0,T))$ and assume that
\begin{equation} \label{withprefactor}
w(t)^\delta \partial_t \big(k\ast [w^\gamma-w_0^\gamma]\big)(t)+\mu(t) w(t)^\eta\le 0,\quad \mbox{a.a.}\,t\in (0,T).
\end{equation}
Then there holds
\begin{equation} \label{withoutprefactor}
\partial_t \big(k\ast [w^\gamma-w_0^\gamma]\big)(t)+\mu(t) w(t)^{\eta-\delta}\le 0,\quad \mbox{a.a.}\,t\in (0,T).
\end{equation}
\end{lemma}
{\em Proof.}
Set $v(t)=w(t)^\gamma$, $v_0=w_0^\gamma$, and define for $\varepsilon>0$ the function $\psi_\varepsilon:[0,\infty)\rightarrow [0,\infty)$ by
\[
\psi_\varepsilon(y)=\int_0^y \frac{r^{\delta/\gamma}}{r^{\delta/\gamma}+\varepsilon}\,\dif r,\quad y\ge 0.
\]
Clearly, $\psi_\varepsilon \in C^1([0,\infty))$ with
\[
\psi_\varepsilon'(y)=\frac{y^{\delta/\gamma}}{y^{\delta/\gamma}+\varepsilon}=1-\frac{\varepsilon}{
y^{\delta/\gamma}+\varepsilon},\quad y\ge 0,
\]
and $\psi_\varepsilon$ is convex. Dividing \eqref{withprefactor} by $w(t)^\delta+\varepsilon$ we get
\[
\psi_\varepsilon'\big(v(t)\big)\, \partial_t \big(k\ast [v-v_0]\big)(t)+\frac{\mu(t) w(t)^\eta}
{w(t)^\delta+\varepsilon}\le 0,
\]
which further yields
\[
\partial_t \big(k\ast [\psi_\varepsilon(v)-\psi_\varepsilon(v_0)]\big)(t)+\frac{\mu(t) w(t)^\eta}
{w(t)^\delta+\varepsilon}\le 0,
\] 
by the convexity inequality \eqref{FIconvex}. 

Next, let $\varphi\in C^\infty_0((0,T))$ be a nonnegative test function. Multiplying the last inequality by 
$\varphi$, integrating over $(0,T)$ and integrating by parts gives
\[
\int_0^T \Big(-(k\ast [\psi_\varepsilon(v)-\psi_\varepsilon(v_0)])(t)\,\partial_t \varphi(t)
+\frac{\mu(t) w(t)^\eta}
{w(t)^\delta+\varepsilon}\,\varphi(t)\Big) \,\dif t\le 0.
\]
Sending $\varepsilon\to 0+$ and using that $\lim_{\varepsilon\to 0+}\psi_\varepsilon(y)= y$ for all $y\ge 0$ 
as well as 
\[
|\psi_\varepsilon(v(t))-\psi_\varepsilon(v_0)|\le \psi_\varepsilon(v(t))+\psi_\varepsilon(v_0)\le v(t)+v_0,\quad t\in [0,T] 
\]
and
\[
\Big|\frac{\mu(t) w(t)^\eta}
{w(t)^\delta+\varepsilon}\Big|\le |\mu(t)| w(t)^{\eta-\delta},\quad  \mbox{a.a.}\,t\in (0,T),
\]
we then obtain by Lebesgue's theorem on dominated convergence that
\[
\int_0^T \Big(-(k\ast [v-v_0])(t)\,\partial_t \varphi(t)
+\mu(t) w(t)^{\eta-\delta}
\,\varphi(t)\Big) \,\dif t\le 0,
\]
and thus
\[
\int_0^T \Big(\partial_t \big(k\ast [v-v_0]\big)(t)
+\mu(t) w(t)^{\eta-\delta}
\Big)\varphi(t) \,\dif t\le 0,
\]
for all nonnegative test functions $\varphi$. The assertion of the lemma follows now by the fundamental lemma of the calculus
of variations.
\hfill $\square$

\subsection{Admissible relative entropies and related inequalities}\label{sec:addmissible_entropies}

Recall that a relative entropy $H(u)$ is induced by a generating function $\phi$, 
cf.~\eqref{relative_entropy}.
In this subsection we describe the class of admissible functions $\phi$ and discuss some important
examples. The following definition has been taken from~\cite{AMTU}. 
\begin{defi}\label{rel_entropy_def}
Let \(I=(0,\infty)\) and \(\phi\in C(\overline{I})\cap C^4(I)\) satisfy the conditions
\begin{align*}
\phi(1)&=0,\\
\phi''\not\equiv 0,\quad\phi''&\ge 0 \quad\text{on $I$},\\
(\phi^{(3)})^2&\le\frac12\phi''\phi^{(4)} \quad\text{on $I$}. 
\end{align*}
Let \(u_1,u_2\in L^1(\iR^d)\) be positive functions with \(\int_{\iR^d} u_1\dif x=\int_{\iR^d} u_2\dif x=1\). Then
\begin{equation}\label{relative_entropy_def}
E_\phi(u_1|u_2)=\int_{\iR^d}\phi\left(\frac{u_1(x)}{u_2(x)}\right)u_2(x)\dif x
\end{equation}
is called an admissible relative entropy of \(u_1\) with respect to \(u_2\) with generating function \(\phi\).
\end{defi}

Comparing this definition with~\eqref{relative_entropy}, we see that the relative entropies 
considered in this paper are of the form
\[
H(u)=H_\phi(u)=E_\phi(u|u_\infty),
\]
where $u_\infty$ is the unique equilibrium of~\eqref{FP_eq} given in \eqref{Gibbs_distribution}.

Important examples of admissible generating functions are the power type function
\begin{equation}\label{entropy_function_p_second}
\phi_\beta(x)=x^\beta-1-\beta(x-1),\quad 1<\beta\le 2,
\end{equation}
and
\begin{equation}\label{logarithmic_entropy}
\phi(x)=x\big(\log(x)-1\big)+1.
\end{equation}
Note that the former takes a simple form in the special case \(\beta=2\). In fact,
\[
\phi_2(x)=x^2-1-2(x-1)=(x-1)^2.
\]
The logarithmic function~\eqref{logarithmic_entropy} leads to the {\em Boltzmann entropy}
\begin{equation}\label{Boltzmann_entropy}
H(u)=\int_{\iR^d} u \log\left(\frac{u}{u_\infty}\right)\,\dif x
\end{equation}
of the probability density $u$ and can be viewed as a 
limiting case of $\phi_\beta$ as $\beta\to 1+$ in the sense that
\[
\lim_{\beta\to 1+}\frac{\phi_{\beta}(x)}{\beta-1}=x\big(\log(x)-1\big)+1.
\]

The following inequality provides a control of the $L^1$ distance of two probability
densities $f,g$ by the relative entropy $E_\phi(f|g)$. It can be found in
\cite[Section 2.2]{AMTU}, see also \cite[Theorem A.3]{jungel2016}.
\begin{satz}[General Csisz\'{a}r-Kullback-Pinsker inequality]\label{CKP_ineq}
Let \(f,g\in L^1(\iR^d)\) be positive functions with unit mass. If \(\phi\) satisfies the assumptions
of Definition~\ref{rel_entropy_def}, then
\[
\|f-g\|_{L^1(\iR^d)}^2\le\frac{2}{\phi''(1)}\int_{\iR^d}\phi\left(\frac{f}{g}\right)g\,\dif x.
\]
\end{satz}

We also need the following version of the convex Sobolev inequality
(see~\cite[Corollary 2.1]{jungel2016}).
\begin{satz} [Convex Sobolev inequality]\label{CSI}
Let $\phi$ be an admissible generating function in the sense of Definition~\ref{rel_entropy_def}
and assume that 
$V\in C^2(\iR^d)$ satisfies the condition 
\eqref{potential_convexity_cond} for some $\lambda>0$ and is such that $e^{-V}\in L^1(\iR^d)$.
Defining $u_\infty$ as in
\eqref{Gibbs_distribution} there holds
\begin{equation}\label{convex_sobolev_ineq}
\int_{\iR^d}\phi\left(\frac{u}{u_\infty}\right)u_\infty\dif x\le\frac{1}{2\lambda}\int_{\iR^d}\phi''\left(\frac{u}{u_\infty}\right)\left|\nabla\frac{u}{u_\infty}\right|^2 u_\infty\,\dif x
\end{equation}
for all nonnegative integrable functions \(u\) for which the integrals are defined.
\end{satz}

\section{Proof of Theorem \ref{thm_main_result}}
\label{ProofMR}

\subsection{Part A - general entropies}

Let the assumptions of Theorem \ref{thm_main_result} (Part A) be satisfied. We set
\begin{equation}\label{FP_oper}
\mathcal{L}^*u:=\Delta u+\nabla\cdot(u \nabla V).
\end{equation}
Observe that
\[
\nabla u_\infty=-u_\infty \nabla V,
\]
and thus
\[
\nabla \cdot \Big( u_\infty \nabla \big(\frac{u}{u_\infty}\big)\Big)=\nabla\cdot 
\Big( \nabla u+u_\infty u \big(-\frac{1}{u^2_\infty}\big)\nabla u_\infty \Big)=\mathcal{L}^*u.
\]

Setting
\[
v:=\frac{u}{u_\infty},\quad v_0:=\frac{u_0}{u_\infty},
\]
equation~\eqref{FP_eq} can be rewritten as
\begin{equation} \label{fp_rewritten}
\partial_t \big(k\ast [v-v_0]\big) u_\infty-\nabla \cdot \big( u_\infty \nabla v\big)=0
\end{equation}
Multiplying \eqref{fp_rewritten} by $\phi'(v)$, integrating over $\iR^d$ and integrating by parts we obtain
\begin{equation}\label{FP_tested2}
\int_{\iR^d}\phi'(v)\partial_t \big(k*[v-v_0]\big) u_\infty\dif x=-\int_{\iR^d}\phi''(v)|\nabla v|^2u_\infty\dif x.
\end{equation}

The right-hand side of~\eqref{FP_tested2} can be estimated from above by means of the convex Sobolev inequality, Theorem \ref{CSI}. This gives
\begin{equation}\label{FP_tested2a}
\int_{\iR^d}\phi'(v)\partial_t \big(k*[v-v_0]\big) u_\infty\dif x\le 
-2\lambda \int_{\iR^d}\phi(v)\,u_\infty\dif x=-2\lambda H(u).
\end{equation}
For the left-hand side of \eqref{FP_tested2a} we use the fundamental identity
in the form of the convexity inequality \eqref{FIconvex} ($\phi$ is convex!), thereby obtaining
\[
\int_{\iR^d} \partial_t \big(k*[\phi(v)-\phi(v_0)]\big) u_\infty\dif x\le 
-2\lambda H(u),
\]
which is equivalent to
\begin{equation} \label{gen_ent_inequ}
\partial_t\big(k\ast \big[H(u)-H(u_0)\big]\big)+2\lambda H(u)\le 0.
\end{equation}
The first assertion from Part A, estimate \eqref{decay_rate_entropyA}, follows now 
from \eqref{gen_ent_inequ} by the comparison principle (cf.~\cite[Section 2.3]{VZ}).

The $L^1$-estimate \eqref{decay_rate_solutionA} is a consequence of \eqref{decay_rate_entropyA}
and the Csisz\'{a}r-Kullback-Pinsker inequality stated in Theorem \ref{CKP_ineq}.
This finishes the proof of Part A of Theorem \ref{thm_main_result}.

\subsection{Part B - power type entropy} \label{SPartB}

The proof of Part B is much more involved. In order to obtain the desired inequality \eqref{decay_rate_entropyB} one also has to exploit the third term in the fundamental identity \eqref{fundamental_identity}. 

The basic idea of the proof is inspired by~\cite[Lemma 3.1]{VZ}. We split the left-hand side of~\eqref{FP_tested2} by writing
\begin{align} \label{I1I2def}
\int_{\iR^d}\phi'(v)\partial_t \big(k*[v-v_0]\big) u_\infty\dif x & = 
\int_{\iR^d}\phi'(v)\partial_t(k*v) u_\infty\dif x - \int_{\iR^d}\phi'(v) v_0 u_\infty\dif x\,k(t)\nonumber \\
& =:I_1(t)-I_2(t).
\end{align}
Applying the fundamental identity~\eqref{fundamental_identity} pointwise in $\iR^d$ 
we have at time $t>0$
\begin{equation}\label{fundamental_identity_for_entropy}
\begin{split}
I_1(t) & = \int_{\iR^d}\partial_t\big(k*\phi(v)\big)u_\infty\dif x+\int_{\iR^d}\big(-\phi(v)+\phi'(v)v\big)u_\infty\dif x\,k(t)\\
&\quad +\int_{\iR^d}\int_0^t\Big(\phi(v(t-s))-\phi(v(t))-\phi'(v(t))[v(t-s)-v(t)]\Big)(-k'(s))\dif s\, u_\infty\dif x. 
\end{split}
\end{equation}
Note that \(v\) depends both on \(t\) and \(x\) although we have used the notation \(v(t)\).
Since~\eqref{fundamental_identity_for_entropy} contains different time levels and the formula
is quite long, it is more convenient to adopt the shorthand notation \(v(t)\) and $u(t)$ instead of \(v(t,x)\)
and $u(t,x)$, respectively.
 
The first term on the right hand side of~\eqref{fundamental_identity_for_entropy} is simply
\[
\partial_t(k*H(u))(t),
\]
in view of the definition of the relative entropy~\eqref{relative_entropy}.

For the second term we use conservation of mass, which gives
\[
\int_{\iR^d}v(t)u_\infty\dif x=\int_{\iR^d}u(t)\dif x=\int_{\iR^d}u_\infty\dif x=1.
\]
This allows to add arbitrary multiples of \(v-1\) inside the spatial integral with weight function $u_\infty$. 
Recalling that
\[
\phi(x)=\phi_\beta(x)=x^\beta-1-\beta(x-1),\quad 1<\beta\le 2,
\]
we thus have at time $t>0$
\[
\begin{split}
\int_{\iR^d}\big(-\phi(v) & +\phi'(v)v\big)u_\infty\dif x\,k(t)  = 
\int_{\iR^d}\big(-\phi(v)+\beta(v^{\beta-1}-1)v\big)u_\infty\dif x\, k(t) \\
&=-H(u(t))k(t)+\beta\int_{\iR^d}\Big(\left(v^\beta-v\right)-(\beta-1)[v-1]\Big)u_\infty\dif x\, k(t)\\
&=(\beta-1)H(u(t))k(t).
\end{split}
\]

The third term on the right-hand side of \eqref{fundamental_identity_for_entropy} is the most  difficult one to handle. First of all,
we have
\[
\begin{split}
& \int_{\iR^d}\int_0^t\Big(\phi(v(t-s))-\phi(v(t))-\phi'(v(t))[v(t-s)-v(t)]\Big)(-k'(s))\dif s\, u_\infty\dif x \\
=&\int_0^t \big(H(u(t-s))-H(u(t))\big)(-k'(s))\dif s\\
&-\beta\int_0^t\int_{\iR^d}\big(v(t)^{\beta-1}-1\big)[v(t-s)-v(t)]u_\infty\dif x\,(-k'(s))\,\dif s.
\end{split}
\]
The first integral on the right hand side is already expressed in terms of the entropy, so this is a good term and we are left with the second integral. Set
\[
I_3(t,s)=\int_{\iR^d}\big(v(t)^{\beta-1}-1\big)[v(t-s)-v(t)]u_\infty\dif x.
\]
Using again conservation of mass, we have
\[
\begin{split}
I_3(t,s) &=\int_{\iR^d}\left(v(t)^{\beta-1}v(t-s)-v(t)^\beta-v(t-s)+v(t)\right)u_\infty\dif x\\
&= -\int_{\iR^d} \big(v(t)^\beta-1-\beta [v(t)-1] \big)u_\infty\,\dif x\\
& \quad +\int_{\iR^d} \Big(-(\beta-1)[v(t)-1]+v(t)^{\beta-1} v(t-s)-v(t-s)\Big)u_\infty\,\dif x\\ 
&=-H(u(t))+\int_{\iR^d}\big(v(t)^{\beta-1}v(t-s)-v(t)\big)u_\infty\,\dif x.
\end{split}
\]
The last integral can be estimated from above by an expression which only involves two entropy terms. Indeed we can prove the following result. 

\begin{lemma}\label{crucial_lemma_entropy}
Let $1<\beta\le 2$. Then for the entropy~\eqref{relative_entropy} with \(\phi:=\phi_\beta\) given by~\eqref{entropy_function_p} there holds
\begin{equation}\label{crucial_inequality_entropy}
\int_{\iR^d}\big(v(t)^{\beta-1} v(t-s)-v(t)\big)u_\infty\dif x\le H(u(t))^{\frac{\beta-1}{\beta}}H(u(t-s))^{\frac{1}{\beta}}
\end{equation}
\end{lemma}

Before entering the proof we make a couple of remarks. Due to conservation of mass we can add arbitrary multiples of 
\((v(t)-1)u_\infty\) and \((v(t-s)-1)u_\infty\) to the integrand on the left-hand side of~\eqref{crucial_inequality_entropy}
without changing the value of the integral. The idea is to use H\"{o}lder's inequality but its direct use is not possible. So we have to modify the left-hand side of~\eqref{crucial_inequality_entropy} appropriately. It
turns out that it is useful to rewrite the left-hand side in the form
\[
\begin{split}
&\int_{\iR^d}\Big(v(t)^{\beta-1}v(t-s)-v(t)+(2-\beta)[v(t)-1]-[v(t-s)-1]\Big)u_\infty\dif x\\
=&\int_{\iR^d}\Big(v(t)^{\beta-1}v(t-s)+(1-\beta)v(t)-v(t-s)+\beta-1)\Big)u_\infty\dif x.
\end{split}
\] 

Then to show~\eqref{crucial_inequality_entropy} we use the following pointwise estimate.

\begin{lemma}\label{crucial_pointwise_lemma}
Let $1<\beta\le 2$. Then for the function \(\phi:=\phi_\beta\) defined by~\eqref{entropy_function_p} there holds
\begin{equation}\label{crucial_pointwise_estimate}
x^{\beta-1}y+(1-\beta)x-y+\beta-1\le\phi(x)^{\frac{\beta-1}{\beta}}\phi(y)^{\frac{1}{\beta}},\quad x,y\ge 0.
\end{equation}
\end{lemma}


{\em Proof.}
We distinguish different cases. We begin with \(\beta=2\). In this case the generating function is simply
\[
\phi(x)=(x-1)^2.
\]
On the other hand, the left-hand side of~\eqref{crucial_pointwise_estimate} is
\[
xy-x-y-1=(x-1)(y-1)\le |x-1|\,|y-1|,
\]
and so~\eqref{crucial_pointwise_estimate} clearly holds.

We consider now the case \(1<\beta<2\).
The left-hand side of~\eqref{crucial_pointwise_estimate} can be written in the form
\begin{equation}\label{lhs_line}
(x^{\beta-1}-1)y+(1-\beta)(x-1),
\end{equation}
which is, for any fixed $x\neq 1$, a first order polynomial in $y$.
If $x=1$, the expression in \eqref{lhs_line} vanishes and the desired inequality is trivially true.
We now distinguish the two cases $x<1$ and $x>1$.

{\em The case \(x<1\).}  In this case~\eqref{lhs_line} is strictly decreasing in $y$, so obviously for large values of \(y\) it is negative and we are done. We have
\[
(x^{\beta-1}-1)y+(1-\beta)(x-1)\le 0\Leftrightarrow y\ge\frac{(\beta-1)(1-x)}{1-x^{\beta-1}}.
\]
The function \(f(x)=x^{\beta-1}-1\) is strictly concave and negative for $x<1$, and so there holds
\[
\frac{(\beta-1)(1-x)}{1-x^{\beta-1}}=(\beta-1)\frac{1-x}{f(1)-f(x)}<\frac{\beta-1}{f'(1)}=1.
\]
Hence it remains to prove the estimate~\eqref{crucial_pointwise_estimate} for \(y<1\). 

Assuming $y<1$ we study carefully the behaviour of the function
\begin{equation}\label{difference_function}
F(x,y)=\phi(x)^{\frac{\beta-1}{\beta}}\phi(y)^{\frac{1}{\beta}}-\big(x^{\beta-1}y+(1-\beta)x-y+\beta-1\big).
\end{equation}
On the diagonal \(y=x\) the function $F$ is zero:
\begin{align*}
F(x,x) & =\phi(x)^{\frac{\beta-1}{\beta}}\phi(x)^{\frac{1}{\beta}}-\big(x^\beta+(1-\beta)x-x+\beta-1\big)\\
& =\phi(x)-(x^{\beta}-1-\beta( x-1))=0.
\end{align*}
We will prove that \(F(x,y)>0\) whenever \(x\neq y\). Let us first calculate the critical points of \(F(x,\cdot)\) with $x<1$
being fixed. We have
\[
\frac{\partial F}{\partial y}=\phi(x)^{\frac{\beta-1}{\beta}}\phi(y)^{\frac{1-\beta}{\beta}}(y^{\beta-1}-1)-(x^{\beta-1}-1),
\]
which is clearly zero, if \(y=x\). To prove that there are no other zeros, it is enough to show that \(F(x,\cdot)\) is 
strictly convex. To this end we calculate the second derivative
\[
\frac{\partial^2 F}{\partial y^2}=\phi(x)^{\frac{\beta-1}{\beta}}\phi(y)^{\frac{1-2\beta}{\beta}}(\beta-1)\big(\phi(y)y^{\beta-2}-(y^{\beta-1}-1)^2\big).
\]
The prefactors are positive, so it remains to study the function
\begin{equation}\label{g_function}
g(y)=\phi(y)y^{\beta-2}-(y^{\beta-1}-1)^2=(\beta-1)y^{\beta-2}+(2-\beta)y^{\beta-1}-1.
\end{equation}
The first derivative of \(g\) is
\begin{equation}\label{g_function_derivative}
g'(y)=(2-\beta)(\beta-1)y^{\beta-3}(y-1),
\end{equation}
which is negative for \(y<1\) and so \(g\) is strictly decreasing. Since \(g(1)=0\), we see that
\[
\frac{\partial^2 F}{\partial y^2}>0,
\]
whence \(F(x,\cdot)\) is strictly convex for \(y< 1\) with $x<1$ being fixed arbitrarily. Therefore \(F\) is zero only on the diagonal and positive outside the diagonal and we are done.

{\em The case \(x> 1\)}. In this case~\eqref{lhs_line} is strictly increasing in $y$, whence the left-hand side of~\eqref{crucial_pointwise_estimate} is nonpositive for all
\[
y\le\frac{(\beta-1)(x-1)}{x^{\beta-1}-1}.
\]
Once again, since \(f(x)=x^{\beta-1}-1\) is strictly concave and $x>1$, we have
\[
\frac{(\beta-1)(x-1)}{x^{\beta-1}-1}>\frac{\beta-1}{f'(1)}=1
\]
and thus it remains to prove~\eqref{crucial_pointwise_estimate} for \(y>1\). In this case \(g\) defined by~\eqref{g_function} is strictly increasing, by \eqref{g_function_derivative}. Since \(g(1)=0\), we have
\[
\frac{\partial^2 F}{\partial y^2}>0,
\]
whenever $y>1$.
Therefore \(F(x,\cdot)\) is again strictly convex for \(y> 1\) with $x>1$ being fixed arbitrarily, which completes the proof.
\hfill$\square$

$\mbox{}$

\noindent Now we can proceed to the proof of Lemma~\ref{crucial_lemma_entropy}.

{\em Proof of Lemma~\ref{crucial_lemma_entropy}.} By Lemma~\ref{crucial_pointwise_lemma} and the remarks prior to it
we may estimate as
\[
\int_{\iR^d}\left(v(t)^{\beta-1} v(t-s)-v(t)\right)u_\infty\dif x\le\int_{\iR^d}\phi(v(t))^{\frac{\beta-1}{\beta}}\phi(v(t-s))^{\frac{1}{\beta}}u_\infty\,\dif x.
\]
The asserted inequality follows then by an application of H\"{o}lder's inequality with the exponents \(p=\frac{\beta}{\beta-1}\) and \(p'={\beta}\) with respect to the measure \(\mu(\dif x)=u_\infty(x)\dif x\).
\hfill$\square$

\begin{remark} \label{generalizedLemma}
Observe that in the proof of Lemma~\ref{crucial_lemma_entropy} we only used the properties that
$u(t)=u(t,\cdot)$, $u(t-s)=u(t-s,\cdot)$ and $u_\infty$ are probability densities on $\iR^d$ with
$u_\infty>0$ everywhere and that $v(t)=u(t)/u_\infty$ as well as $v(t-s)=u(t-s)/u_\infty$. Thus the 
argument yields a more general result. For example, it shows that if $f_1, f_2, g$ are probability densities on $\iR^d$ with $g>0$ on $\iR^d$ then we have with $h_i:=f_i/g$, $i=1,2$,
\[
\int_{\iR^d} h_1^{\beta-1} h_2 g\,\dif x-1 = \int_{\iR^d} \big( h_1^{\beta-1} h_2- h_1\big) g\,\dif x \le 
E_{\phi_\beta}(f_1|g)^{\frac{\beta-1}{\beta}} E_{\phi_\beta}(f_2|g)^{\frac{1}{\beta}},
\]
where $\phi_\beta$ is as in Lemma~\ref{crucial_lemma_entropy}, cf.\ \eqref{relative_entropy_def} for the definition of the relative entropy. Evidently, as the above proof shows, this statement can be 
generalized further to more general probability spaces.
\end{remark}

$\mbox{}$

\noindent Having Lemma~\ref{crucial_lemma_entropy} at our disposal we can estimate the term $I_3(t,s)$
from above by entropy terms as follows:
\[
I_3(t,s)\le -H(u(t))+H(u(t))^{\frac{\beta-1}{\beta}}H(u(t-s))^{\frac{1}{\beta}}.
\]
Combining this estimate and the above calculations concerning the term $I_1(t)$ we obtain for $t>0$
\begin{align*}
I_1(t) & = \partial_t\big(k\ast H(u)\big)(t)+(\beta-1) H(u(t))k(t)\\
& \quad +\int_0^t \big(H(u(t-s))-H(u(t))\big)(-k'(s))\,\dif s
-\beta\int_0^t I_3(t,s)\,(-k'(s))\,\dif s\\
& \ge \partial_t\big(k\ast H(u)\big)(t)+(\beta-1) H(u(t))k(t)\\
& \quad +\int_0^t \left(H(u(t-s))+(\beta-1)H(u(t))-\beta H(u(t))^{\frac{\beta-1}{\beta}}H(u(t-s))^{\frac{1}{\beta}}\right)(-k'(s))\,\dif s\\
& =: I_4(t).
\end{align*}
 
A key idea is now to apply the fundamental
identity a second time (!) to show that 
\begin{equation} \label{I4fromula}
I_4(t)=\beta H(u(t))^{\frac{\beta-1}{\beta}}\partial_t \big(k* H(u)^{\frac{1}{\beta}}\big)(t).
\end{equation} 

In fact, setting 
\[
w(t)=H(u(t))^{\frac{1}{\beta}}
\]
the fundamental identity~\eqref{fundamental_identity} with $\psi(y)=y^\beta$
gives
\begin{align*}
& \beta  H(u(t))^{\frac{\beta-1}{\beta}} \partial_t   \big(k* H(u)^{\frac{1}{\beta}}\big)(t)
 = \beta w(t)^{\beta-1}\partial_t (k* w)(t)\\
& = \,\partial_t (k\ast w^\beta)(t)+(\beta-1)w(t)^\beta k(t)\\
& \quad\, +\int_0^t\Big(w(t-s)^\beta-w(t)^\beta-\beta w(t)^{\beta-1}[w(t-s)-w(t)]\Big)(-k'(s))\,\dif s\\
& = \,\partial_t \big(k\ast H(u)\big)(t)+(\beta-1) H(u(t)) k(t)\\
& \quad \,+\int_0^t\Big(H(u(t-s))-H(u(t))-\beta H(u(t))^{\frac{\beta-1}{\beta}}
\big[H(u(t-s))^{\frac{1}{\beta}}-H(u(t))^{\frac{1}{\beta}}\big]\Big)(-k'(s))\,\dif s\\
& = I_4(t). 
\end{align*}
Combining \eqref{I4fromula} and $I_1(t)\ge I_4(t)$ we obtain
\begin{equation} \label{I1estimate}
I_1(t)\ge \beta H(u(t))^{\frac{\beta-1}{\beta}}\partial_t \big(k* H(u)^{\frac{1}{\beta}}\big)(t).
\end{equation}

It remains to estimate the term $I_2(t)$, which contains the initial datum. Using again conservation
of mass we have
\begin{align*}
I_2(t) = \int_{\iR^d} \beta\big( v(t)^{\beta-1}-1\big) v_0 u_\infty\,\dif x\,k(t)
= \beta\int_{\iR^d}\left(v(t)^{\beta-1}v_0-v(t)\right)u_\infty\,\dif x\,k(t).
\end{align*}
By virtue of Remark \ref{generalizedLemma} it follows that
\begin{equation} \label{I2estimate}
I_2(t)\le \beta k(t) H(u(t))^{\frac{\beta-1}{\beta}} H(u_0)^{\frac{1}{\beta}}. 
\end{equation}

Combining \eqref{FP_tested2a}, \eqref{I1I2def}, \eqref{I1estimate}, \eqref{I2estimate} we obtain
\begin{equation} \label{entestprefactor}
\beta H(u(t))^{\frac{\beta-1}{\beta}}\partial_t\big(k*(H(u)^{1/\beta}-H(u_0)^{1/\beta})\big)(t)\le -2\lambda H(u(t)),
\quad t>0.
\end{equation}
Invoking Lemma \ref{dividingprefactor} it follows from \eqref{entestprefactor} that
\begin{equation}\label{decay_rate_entropy_final}
\partial_t(k*(H(u)^{1/\beta}-H(u_0)^{1/\beta}))(t)\le -\frac{2\lambda}{\beta} H(u(t))^{1/\beta},
\end{equation}
which in turn implies the asserted estimate \eqref{decay_rate_entropyB},
by the comparison principle (cf.~\cite[Section 2.3]{VZ}).
The second estimate in Part B, inequality \eqref{decay_rate_solutionB}, follows directly from \eqref{decay_rate_entropyB}
by the Csisz\'{a}r-Kullback-Pinsker inequality stated in Theorem \ref{CKP_ineq} with $\phi=\phi_\beta$. 
This finishes the proof of Part B of Theorem \ref{thm_main_result}.

\section{Examples and proof of Corollary \ref{main_result_timefrac}}
In this section we consider several important examples of kernels of type $\mathcal{PC}$
and look at the long-time behaviour of the corresponding relaxation functions $s_\mu$. 
For more details and further examples we refer to \cite[Section 6]{VZ}. 

\begin{bei} \label{beifrac}
The time-fractional case. {\em We consider the pair
\[
(k,l)=(g_{1-\alpha},g_\alpha),
\]
where \(\alpha\in(0,1)\) and \(g_\beta\) is defined by~\eqref{standard_kernel}. Recall that the Laplace transform of $g_\beta$, $\beta>0$, is given by $\widehat{g_\beta}(z)=z^{-\beta}$, Re$\,z>0$, so it is easy to see that $g_{\beta_1}\ast
g_{\beta_2}=g_{\beta_1+\beta_2}$ for all $\beta_1,\beta_2>0$. This implies
\[
(k*l)(t)=(g_{1-\alpha}*g_\alpha)(t)=g_1(t)=1,\quad t>0,
\]
and so it is clear that $(g_{1-\alpha},g_\alpha)\in \mathcal{PC}$. We further have
\[
(1\ast l)(t)=(1\ast g_\alpha)(t)=g_{1+\alpha}(t)=\frac{t^\alpha}{\Gamma(1+\alpha)},
\]
which together with~\eqref{smubounds} shows that for $\mu\ge 0$ the corresponding relaxation
function $s_\mu$ satisfies
\begin{equation} \label{timefracbounds}
\frac{1}{1+\mu \Gamma(1-\alpha) t^\alpha}\le s_\mu(t)\le 
\frac{1}{1+\mu [\Gamma(1+\alpha)]^{-1} t^\alpha},\quad t\ge 0.
\end{equation}
Corollary \ref{main_result_timefrac} now follows from Theorem \ref{thm_main_result} and \eqref{timefracbounds}.  

We remark that
\[
s_\mu(t)=E_\alpha(-\mu t^\alpha),\quad\mbox{where}\;E_\alpha(z):=\sum_{j=0}^\infty
\,\frac{z^j}{\Gamma(\alpha j+1)}\,,\;z\in \iC,
\]
is the well-known Mittag-Leffler function (see e.g.\ \cite{KST}).

If one replaces the above pair of kernels by
\[
k(t)=g_{1-\alpha}(t)e^{-\gamma t},\quad l(t)=g_{\alpha}(t)e^{-\gamma
t}+\gamma(1\ast[g_{\alpha}e^{-\gamma\cdot}])(t), \quad t>0,
\]
with $\alpha\in (0,1)$ and $\gamma>0$, then again $(k,l)\in \mathcal{PC}$ and
the associated relaxation function decays exponentially, see \cite[Example 6.2]{VZ}.
}
\end{bei}
\begin{bei} \label{multiterm}
Multiterm fractional derivative. {\em We consider the sum
of two fractional derivatives. Let $0<\alpha<\beta<1$ and
\[
k(t)=g_{1-\alpha}(t)+g_{1-\beta}(t), \quad t>0.
\]
Clearly, $k$ is completely monotone and $k(0+)=\infty$. It follows by Theorem 5.4 in Chapter 5 of \cite{GLS} that the kernel
$k$ has a resolvent $l\in L^1_{loc}(\iR_+)$ of the first kind,
that is $k\ast l=1$ on $(0,\infty)$, and this resolvent is
completely monotone as well. In particular $(k,l)\in {\cal PC}$. Since
\[
\widehat{1\ast l}\,(z)=\,\frac{1}{z}\,\frac{1}{z^\alpha+z^\beta}\,\sim \frac{1}{z^{1+\alpha}}\quad \mbox{as}\;z\to 0,
\]
the Karamata-Feller Tauberian theorem (see \cite{Feller}) implies that
$(1\ast l)(t)\sim g_{1+\alpha}(t)$ as $t\to \infty$. 
On the other hand, $k(t)\sim g_{1-\alpha}(t)$ as $t\to \infty$, and so, using 
\eqref{smubounds}, we see that the relaxation function $s_\mu(t)$ has the same algebraic
decay (as $t\to \infty$) as in Example \ref{beifrac}; the decay rate is determined
by the fractional derivative of lower order. 

These considerations extend trivally to kernels $k(t)=\sum_{j=1}^m \delta _j g_{1-\alpha_j}(t)$ with $\delta_j>0$ and
$0<\alpha_1<\alpha_2<\ldots<\alpha_m<1$.
}
\end{bei}
\begin{bei} \label{beislow}
The distributed order case (an example of ultraslow diffusion). {\em We consider the pair
\[
k(t)=\int_0^1 g_\beta(t)\,\dif\beta,\quad l(t)=\int_0^\infty \,\frac{e^{-st}}{1+s}\,\dif s,\quad t>0.
\]
Both kernels are nonnegative and nonincreasing, and there holds (see \cite[Example 6.5]{VZ})
\[
\hat{k}(z)=
\,\frac{z-1}{z\log z},\quad  \hat{l}(z)=\,\frac{\log z}{z-1}\,,\quad  \mbox{Re}\,z>0.
\]
Thus $(k,l)\in \mathcal{PC}$. There exists a number $T_1>1$ such that
\[
\frac{1}{2k(t)}\,\le \log t\le 2(1\ast l)(t),\quad t\ge T_1,
\]
see \cite[Example 6.5]{VZ}. In view of \eqref{smubounds} this shows that for $\mu>0$ the relaxation function $s_\mu(t)$ can be estimated for large times from above and from below by terms of the form
$c/\log  t$ with some constant $c>0$. So we get a logarithmic decay, which is certainly slower than
in the time-fractional case for any $\alpha\in (0,1)$.  
}
\end{bei}
\section{Optimality of the entropy decay rate in the case $\beta=2$} \label{optHilbert}
We consider the situation of Part B in Theorem \ref{thm_main_result} with $\beta=2$.
The goal of this section is to show that in this case the estimate
\eqref{decay_rate_entropyB} is in general the best possible one can have.
Since
\[
H(u)=\int_{\iR^d} \left(\frac{u}{u_\infty}-1\right)^2 u_\infty\,\dif x
=\int_{\iR^d} \big(u-u_\infty)^2 u_\infty^{-1}\,\dif x
\]
and $u_\infty\in L^1(\iR^d)$, the assumption $H(u_0)<\infty$ means exactly that
$u_0$ belongs to the weighted $L^2$-space
\[
L^2(u_\infty^{-1}):=L^2(\iR^d; u_\infty^{-1})=\{w:\,\iR^d\rightarrow \iR\;\mbox{measurable s.t.}\;
||w||^2_{L^2(u_\infty^{-1})}:=\int_{\iR^d} w^2 u_\infty^{-1}\, \dif x<\infty\}.
\]
The entropy decay estimate \eqref{decay_rate_entropyB} can be written as
\begin{equation} \label{decayestimateHilbert}
||u(t)-u_\infty||_{L^2(u_\infty^{-1})}\le s_\lambda(t)||u_0-u_\infty||_{L^2(u_\infty^{-1})},\quad t>0.
\end{equation}
In terms of $v:=u/u_\infty$ and $v_0:=u_0/u_\infty$ the estimate \eqref{decay_rate_entropyB}
takes the form
\begin{equation} \label{decayestimateHilbert2}
||v(t)-1||_{L^2(u_\infty)}\le s_\lambda(t)||v_0-1||_{L^2(u_\infty)},\quad t>0,
\end{equation}
with the weighted $L^2$-space $L^2(u_\infty):=L^2(\iR^d; u_\infty)$.
The convex Sobolev inequality \eqref{CSI} becomes the Poincar\'e inequality w.r.t.\ the measure
$d\mu=u_\infty dx$ for the
function $v-1$ (which has zero mean); more precisely,
\begin{equation} \label{poincarev}
\int_{\iR^d} (v-1)^2 u_\infty\,\dif x\le \frac{1}{\lambda}\,\int_{\iR^d} |\nabla v|^2 u_\infty\,\dif x.
\end{equation}

Anagolously to the classical case, which is discussed e.g.\ in \cite{Pavliotis}, $u$ solves the
non-local Fokker-Planck equation~\eqref{FP_eq} with initial value $u_0$ if and only if
$v=u/u_\infty$ solves the non-local {\em backward Kolmogorov equation}
\begin{equation} \label{Kolmback}
\partial_t\big(k*[v-v_0]\big)-\mathcal{L}v=0,\quad t>0,\,x\in \iR^d,
\end{equation}
with initial value $v_0=u_0/u_\infty$, where the operator $\mathcal{L}$ is given by
\[
\mathcal{L}f=\Delta f-\nabla V\cdot\nabla f.
\]
It is known that the operator $\mathcal{L}$ is the $L^2(\iR^d)$ adjoint of the Fokker-Planck
operator $\mathcal{L}^*$ appearing in \eqref{FP_eq} (see \eqref{FP_oper} for its definition).
Moreover, $\mathcal{L}$
is self-adjoint with respect to the weighted \(L^2\) inner product
\[
(f,g)_{L^2(u_\infty)}=\int_{\iR^d}f g u_\infty\,\dif x.
\] 
The kernel of $\mathcal{L}$ consists of constants and
\[
(\mathcal{L} f,f)_{L^2(u_\infty)}=-\int_{\iR^d} |\nabla f |^2 u_\infty\,\dif x,
\]
for all sufficiently smooth functions $f$. The spectrum of $-\mathcal{L}$ is discrete consisting
of eigenvalues
\[
0=\lambda_0<\lambda_1<\lambda_2<\dots,
\]
and an orthonormal basis of \(L^2(u_\infty)\) can be built of corresponding
(normalized) eigenfunctions \(\{\varphi_k\}_{k=0}^\infty\) with \(\varphi_0\equiv 1\), 
see~\cite[Chapter 4]{Pavliotis}.
For example, in the case $V(x)=\frac{1}{2}\,|x|^2$, $\mathcal{L}$ becomes the well known
Ornstein-Uhlenbeck operator
\begin{equation} \label{OUop}
\mathcal{L}f(x)=\Delta f(x)-x\cdot \nabla f(x), 
\end{equation}
the spectrum of $-\mathcal{L}$ coincides with $\iN_0$, and there is an orthonormal basis of \(L^2(u_\infty)\) consisting
of $d$-dimensional normalized Hermite polynomials which are eigenfunctions of $-\mathcal{L}$,
see \cite[Chapter 9]{LB}.

Let us first consider the case of arbitrary initial values $u_0$ with finite entropy. Then we may choose as initial value for $u$ 
\[
u_0(x)=\big(1+\varphi_1(x)\big)u_\infty(x),
\]
since
\[
H(u_0)=\int_{\iR^d} \varphi_1^2 u_\infty\,\dif x<\infty,
\]
due to $\varphi_1\in L^2(u_\infty)$. We then have $v_0=1+\varphi_1$, and the corresponding
solution $v$ of the backward Kolmogorov equation \eqref{Kolmback} is given by
\begin{equation}\label{backward_solution}
v(t,x)=s_{\lambda_1}(t)\varphi_1(x)+1,
\end{equation}
since
\begin{align*}
\partial_t\big(k\ast [v-v_0]\big)(t,x) & = \varphi_1(x) \partial_t\big(k\ast [s_{\lambda_1}-1]\big)(t,x)\\
& = -\varphi_1(x)\lambda_1 s_{\lambda_1}(t) \\
& = \mathcal{L}\varphi_1(x) s_{\lambda_1}(t) =  \mathcal{L}v(t,x).
\end{align*}
For this solution, there holds
\[
||v(t)-1||_{L^2(u_\infty)}= s_{\lambda_1}(t)||v_0-1||_{L^2(u_\infty)},\quad t\ge 0,
\]
that is we have equality in \eqref{decayestimateHilbert2}.

In the above example, the initial value $u_0$ and hence the solution \(u\) of~\eqref{FP_eq} can 
assume negative values as can be seen e.g. in the case of the one-dimensional Ornstein-Uhlenbeck operator (see~\cite[Section 4.4]{Pavliotis}).

However, we may easily modify the example presented above to cover also the case of positive solutions. Indeed, take a probability density \(u_0\) such that \(H(u_0)<\infty\). Then \(v_0\) is an element of \(L^2(u_\infty)\) and can be expanded as
\[
v_0(x)=1+\sum_{k=1}^\infty c_k\varphi_k(x).
\]
Then the solution for the backward Kolmogorov equation~\eqref{Kolmback} is given by the series
\[
v(t,x)=1+\sum_{k=1}^\infty c_k s_{\lambda_k}(t)\varphi_k(x).
\]
If \(u_0\) is chosen such that \(c_1\neq 0\), then by Parseval's identity
\[
\|v(t)-1\|_{L^2(u_\infty)}^2=\sum_{k=1}^\infty c_k^2\big(s_{\lambda_k}(t)\big)^2\ge c_1^2 \big(s_{\lambda_1}(t)\big)^2,
\]
which shows that in general the entropy cannot decay faster than \((s_{\lambda_1}(t))^2\) (up to some
positive constant).

In the case of the Ornstein-Uhlenbeck operator \eqref{OUop} in $\iR^d$ the potential
$V(x)=\frac{1}{2}\,|x|^2$ satisfies $\nabla^2 V(x)=Id$ and thus the 
Bakry-Emery condition \eqref{potential_convexity_cond} holds with $\lambda=1$. 
Hence $\lambda=\lambda_1=1$. 

We finally point out that it is an open question whether the obtained decay rates in Theorem \ref{thm_main_result}
are optimal if $\phi\neq \phi_2$.

$\mbox{}$




\section{Time-discrete Fokker-Planck equation}

The purpose of this section is to show that our method also applies to the time-discretized 
Fokker-Planck equation
\begin{equation}\label{FP_discretized}
\frac{1}{\tau}\big(u(t_n,x)-u(t_{n-1},x)\big)-(\mathcal{L}^*u)(t_n,x)=0,\quad n\in \iN,\,x\in \iR^d,
\end{equation}
with initial condition
\begin{equation} \label{discrete_initial_cond}
u(0,x)=u_0(x),\quad x\in \iR^d.
\end{equation}
Here we use the {\em backward difference scheme} for the usual time derivative \(\partial_t\) with
time step \(\tau\), that is \(t_n=n\tau\), $n\in \iN_0$. The Fokker-Planck operator
$\mathcal{L}^*$ is as before, see \eqref{FP_oper} for its definition. 

By means of the crucial estimate from Remark \ref{generalizedLemma}
for the generating function $\phi_\beta$ and the associated entropy
we are able to improve, in the case of a power type entropy, the rate of (exponential) convergence known in the literature for general entropies. So also in the time-discrete
case the decay rates become better with higher
integrability of the initial datum. 

We will prove the following theorem.
\begin{theorem} \label{decaydiscretetime}
Suppose that $V\in C^2(\iR^d)$ satisfies condition 
\eqref{potential_convexity_cond} for some $\lambda>0$ and is such that $e^{-V}\in L^1(\iR^d)$.
Let $u_0\in L^1(\iR^d)$ be a probability density and $u_\infty$ be defined as in
\eqref{Gibbs_distribution}.
Let further $\tau>0$ and $t_n=n \tau$, $n\in \iN_0$. Let $\beta\in (1,2]$ and $\phi_\beta$ be defined as in 
\eqref{entropy_function_p}. Let $H(u)$ be the corresponding relative entropy associated with $u_\infty$.  Assume that $H(u_0)<\infty$ and that the positive solution $u$ of  \eqref{FP_discretized}, \eqref{discrete_initial_cond}
is sufficiently smooth. Then
\begin{equation} \label{discrete_decay_rate_entropyB}
H\big(u(t_n)\big)\le \left(\frac{1}{1+\frac{2\tau \lambda}{\beta}}\right)^{\beta n} H(u_0),\quad n\in \iN.
\end{equation}
Moreover,
\begin{equation} \label{discrete_decay_rate_solutionB}
||u(t_n)-u_\infty||_{L^1(\iR^d)}\le \sqrt{\frac{2}{\beta(\beta-1)}}\left(\frac{1}{1+\frac{2\tau \lambda}{\beta}}\right)
^{\frac{\beta n}{2}}
\sqrt{H(u_0)},\quad n\in \iN.
\end{equation}
\end{theorem}

To compare this result with what is known in the literature, we remark that
if one replaces $\phi_\beta$ by an arbitrary admissible generating function $\phi$, then
it is already known that
\begin{equation} \label{discretegeneralphi}
H\big(u(t_n)\big)\le\left(\frac{1}{1+2\tau\lambda}\right)^n H(u_0),\quad n\in \iN,
\end{equation}  
see~\cite[Remark 5.9]{jungel2016}; we also refer to \cite{AU}. Observe that for $\beta\in (1,2]$
\begin{equation} \label{smallernumber}
\left(\frac{1}{1+\frac{2\tau \lambda}{\beta}}\right)^{\beta }<\frac{1}{1+2\tau\lambda}.
\end{equation}
In fact, setting $\delta=\frac{1}{\beta}$ the function $g(y)=y^\delta$ is strictly concave in $[0,\infty)$ and thus
we have for $a>0$
\[
(1+a)^\delta <1^\delta+\delta \,1^{\delta-1} a=1+\frac{a}{\beta}.
\]
The claim \eqref{smallernumber} then follows easily from this inequality with $a=2\tau \lambda$.

\subsection{The fundamental identity for the discrete time derivative}

For fixed $\tau>0$ we
define the operator $D$ by
\[
Du(t)=\frac{1}{\tau}\,\big(u(t)-u(t-\tau)\big).
\]

Let $I\subset \iR$ be an interval, $\psi \in C^1(I)$, $N\in \iN$, $J:=\{n\tau:\,n=0,\ldots,N\}$,
and $u:\,J\rightarrow I$. Then there holds an analogue of the fundamental identity \eqref{fundamental_identity},
more precisely, we have for all $t\in J\setminus\{0\}$
\begin{equation} \label{discreteFI}
\psi'\big(u(t)\big) (Du)(t)= D\big(\psi(u)\big)(t)+\frac{1}{\tau}\,\Big(\psi\big(u(t-\tau)\big)-\psi\big(u(t)\big)-
\psi'\big(u(t)\big)\big[u(t-\tau)-u(t)\big]\Big).
\end{equation}
This is folklore and follows directly from the definition of $D$. We remark that corresponding identities also hold for 
operators which are discrete in space, e.g.\ the Laplacian on graphs (\cite{DKZ}).

As a direct consequence of \eqref{discreteFI} we obtain for convex functions $\psi$ that
\begin{equation} \label{discreteFIconvex}
\psi'\big(u(t)\big) (Du)(t)\ge  D\big(\psi(u)\big)(t).
\end{equation}

\subsection{Proof of Theorem \ref{decaydiscretetime}}

Let $u$ be the solution of problem \eqref{FP_discretized}, \eqref{discrete_initial_cond}. We set again $v=u/u_\infty$
and $v_0=u_0/u_\infty$. Multiplying~\eqref{FP_discretized} by $\phi'(v(t_n))$, integrating over $\iR^d$ and integrating
by parts we obtain (cf.\ the beginning of Section \ref{ProofMR})
\begin{equation}\label{FP_discretized2}
\tau^{-1} \int_{\iR^d}(v(t_n)-v(t_{n-1}))\phi'(v(t_n))u_\infty\dif x=-\int_{\iR^d}\phi''(v(t_n))|\nabla v(t_n)|^2u_\infty\dif x.
\end{equation}

Let us first recall how the general estimate \eqref{discretegeneralphi} can be derived. Since \(\phi\) is convex, we can apply the convexity inequality \eqref{discreteFIconvex} for the discrete time
derivative to estimate the left-hand side of \eqref{FP_discretized2} as follows
(see also~\cite[Section 5.3]{jungel2016}):
\begin{align}
\int_{\iR^d}(v(t_n)-v(t_{n-1}))\phi'(v(t_n))u_\infty\dif x & \ge\int_{\iR^d}\big(\phi(v(t_n))-\phi(v(t_{n-1}))\big)u_\infty\dif x \nonumber\\
& =H(u(t_n))-H(u(t_{n-1})).  \label{FP_discretized3}
\end{align}

As in the continuous case, the right-hand side of~\eqref{FP_discretized2} can be estimated from above by means of the convex Sobolev inequality, Theorem \ref{CSI}. This yields
\[
\tau^{-1}\big(H(u(t_n))-H(u(t_{n-1}))\big)\le -2\lambda H(u(t_n)),
\]
which shows that $H(u)$ is nonincreasing and which is equivalent to
\[
H(u(t_n))\le\frac{1}{1+2\tau\lambda} H(u(t_{n-1})),\quad n\in \iN.
\]
This inequality in turn implies the general estimate
\eqref{discretegeneralphi}. 

In order to obtain the better estimate for $\phi=\phi_\beta$ we proceed analogously to Subsection \ref{SPartB}.
The key idea is again to use the fundamental identity (now in the form \eqref{discreteFI}) in its full strength.
Denoting again the discrete time derivative by $D$, by \eqref{discreteFI}
we can write the left-hand side of~\eqref{FP_discretized2} in the form
\begin{equation}\label{FP_discretized5}
\begin{split}
\int_{\iR^d}(Dv)(t_n) & \phi'(v(t_n))u_\infty\dif x
= \,\int_{\iR^d} D(\phi(v))(t_n)  u_\infty\dif x\\
&\,+\tau^{-1}\int_{\iR^d}\Big(\phi(v(t_{n-1}))-\phi(v(t_n))-\phi'(v(t_n))\big[v(t_{n-1})-v(t_n)\big]\Big)u_\infty\dif x\\
& =: A_1(n)+A_2(n). 
\end{split}
\end{equation}
Evidently,
\begin{equation} \label{A1term}
A_1(n)=D \big(H(u)\big)(t_n).
\end{equation}
The structure of the second term in \eqref{FP_discretized5} is the same as that of the third term on the right-hand side
of \eqref{fundamental_identity_for_entropy}. So we can follow the line of arguments given in Subsection \ref{SPartB}.
We have
\begin{align*}
A_2(n) & = \tau^{-1}\big(H(u(t_{n-1}))-H(u(t_n))\big)-\frac{\beta}{\tau}\, \int_{\iR^d} \big(v(t_n)^{\beta-1}-1\big)
[v(t_{n-1})-v(t_n)]u_\infty \,\dif x\\
& = \tau^{-1}\big(H(u(t_{n-1}))-H(u(t_n))\big)+\frac{\beta}{\tau}\,H(u(t_n))\\
& \quad\,-\frac{\beta}{\tau}\, \int_{\iR^d} \big(v(t_n)^{\beta-1}
v(t_{n-1})-v(t_n)\big)u_\infty \,\dif x\\
& \ge \tau^{-1}\big(H(u(t_{n-1}))+(\beta-1)H(u(t_n))\big)-\frac{\beta}{\tau}\,H(u(t_n))^{\frac{\beta-1}{\beta}}
H(u(t_{n-1}))^{\frac{1}{\beta}},
\end{align*}
where we used the inequality from Remark \ref{generalizedLemma}. Combining this, \eqref{A1term} and \eqref{FP_discretized5} and setting $w(t)=H(u(t))^{\frac{1}{\beta}}$
we see that
\begin{align*}
\int_{\iR^d} &(Dv)(t_n) \phi'(v(t_n))u_\infty\dif x
  \ge D \big(H(u)\big)(t_n)\\
 & \quad\,+\tau^{-1}\big(H(u(t_{n-1}))+(\beta-1)H(u(t_n))\big)-\frac{\beta}{\tau}\,H(u(t_n))^{\frac{\beta-1}{\beta}}
H(u(t_{n-1}))^{\frac{1}{\beta}}\\
& = D(w^\beta)(t_n)
+\tau^{-1}\Big(w(t_{n-1})^\beta-w(t_n)^\beta-\beta w(t_n)^{\beta-1}\big[w(t_{n-1})-w(t_n)\big]\Big)\\
& = \beta w(t_n)^{\beta-1} (Dw)(t_n),
\end{align*} 
where in the last step we applied the fundamental identity \eqref{discreteFI}.  
   
Altogether, we obtain
\begin{equation} \label{almostend}
\beta w(t_n)^{\beta-1} (Dw)(t_n)\le -2\lambda w(t_n)^\beta,\quad n\in \iN.
\end{equation}
The sequence $(H(u(t_n)))_{n\in \iN_0}$ is nonnegative and nonincreasing, so $(w(t_n))_{n\in \iN_0}$
enjoys the same properties. Consequently, if $w(t_{n_*})=0$ for some $n_*\in \iN_0$ then $w(t_n)=0$ for all
$n\ge n_*$ and thus the asserted entropy estimate trivially holds for all $n\ge n_*$. It remains to look at those
times $t_n$, where $w(t_n)>0$. 

So assume that $w(t_n)>0$ for all $n\in\{0,\ldots,N\}$. Then \eqref{almostend} implies
that
\[
\tau^{-1} \big(w(t_n)-w(t_{n-1})\big)\le -\frac{2\lambda}{\beta}\,w(t_n),\quad n\in \{1,\ldots,N\},
\]  
and thus
\begin{equation}\label{FP_discretized6}
H(u(t_n))\le\left(\frac{1}{1+\frac{2\tau \lambda}{\beta}}\right)^{\beta n} H(u_0),\quad n\in \{1,\ldots,N\}.
\end{equation}
Since $H(u(t_n))=0$ if and only $w(t_n)=0$ and by the remarks following \eqref{almostend}, the last inequality even holds for all $n\in \iN$, thus proving \eqref{discrete_decay_rate_entropyB}. The second assertion of the theorem
follows from \eqref{discrete_decay_rate_entropyB} and the Csisz\'{a}r-Kullback-Pinsker inequality stated in Theorem \ref{CKP_ineq}.



$\mbox{}$

\noindent {\footnotesize {\bf Jukka Kemppainen}, Applied and Computational Mathematics,
Pentti Kaiteran katu 1, PO Box 8000, FI-90014 University of Oulu,
Finland, e-mail: Jukka.T.Kemppainen@oulu.fi





$\mbox{}$

\noindent {\bf Rico Zacher}, Ulm University, Institute of Applied Analysis, 89069 Ulm, Germany,
e-mail: rico.zacher@uni-ulm.de

}

\end{document}